\documentclass[autocontact]{gaceta}

\usepackage[utf8]{inputenc}
\usepackage{etoolbox}
\AtBeginEnvironment{quote}{\itshape}
\usepackage[spanish]{babel}

\journame{La Gaceta de la RSME}
\yearofpublication{0000}
\volume{00}
\issuenumber{0}

\belongstopart{Art\'{\i}culos} 

\title{Zoel García de Galdeano y su presencia en sociedades matemáticas extranjeras de su época}
\author{Pedro J. Miana y Antonio M. Oller-Marcén}
\shorttitle{García de Galdeano en sociedades matemáticas extranjeras de su época}

\contact{Pedro J. Miana, Dpto. de Matem\'aticas-IUMA, Universidad de Zaragoza}
{pjmiana@unizar.es}{http://www.unizar.es/matematicas/personales/pjmiana.htm}
\contact{Antonio M. Oller-Marcén, Centro Universitario de la Defensa de Zaragoza-IUMA}
{oller@unizar.es}{http://cud.unizar.es/aollerm}

\begin{document}

\maketitle

\section{Introducción}

Durante la segunda mitad del siglo XIX y los inicios del siglo XX comenzaron a establecerse por todo el mundo las primeras sociedades matemáticas \cite{FUR2}. Entre ellas, la \textit{London Mathematical Society} en 1865, la \textit{Societé Mathématique de France} en 1872, el \textit{Circolo Matematico di Palermo} en 1884, la \textit{Deutsche Mathematiker-Vereinigung} en 1890, la \textit{American Mathematical Society} en 1894, la \textit{Sociedad Matemática Española} en 1911 o la \textit{Mathematical Association of America} en 1915.

Este fenómeno, junto con la proliferación de revistas matemáticas que en muchos casos estaban asociadas a las sociedades anteriores \cite{GIS}, la convocatoria de los primeros congresos internacionales de matemáticas \cite{ALB} o la creación del CIEM (actualmente conocido como ICMI) \cite{HOD} ponen de manifiesto la creciente y progresiva creación de una cierta conciencia de clase respecto de la profesión de matemático \cite{PAR} y de la consecuente necesidad de comunicación a un nivel internacional.

En España, pese al cierto retraso existente en cuanto a las matemáticas enseñadas y aprendidas en la universidad \cite{HOR3}, surgieron algunos esfuerzos hacia la introducción de nuevas ideas. Una de las figuras más destacadas a este respecto, si no la más destacada, fue la de D. Zoel García de Galdeano (Figura 1). 

Zoel García de Galdeano nació en Pamplona en 1846 y falleció en Zaragoza en 1924. Estudió Ciencias Exactas en la Universidad de Zaragoza y desde 1876 hasta 1889 recorrió España como profesor de la Enseñanza Secundaria en Escuelas e Institutos hasta que consiguió la Cátedra de Geometría Analítica en la Universidad de Zaragoza donde permaneció durante casi treinta años hasta su jubilación en 1918.

Además de un divulgador incansable, de un prolífico autor de libros y folletos \cite{HOR2} o de editar la primera revista matemática española (Hormigón, 1981), García de Galdeano fue el primer matemático español con una clara proyección internacional y una nutrida red de relaciones con importantes matemáticos de Europa y América.

\begin{figure}[h]
\center
\includegraphics{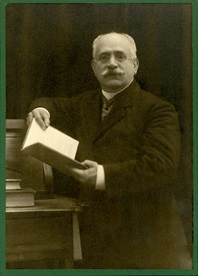}
\caption{Zoel García de Galdeano hacia 1900, retrato de José Yanguas. Fuente: Biblioteca del Dep. de Matemáticas de la Universidad de Zaragoza.}
\end{figure}

Esta vocación internacional de García de Galdeano ha sido puesta de manifiesto habitualmente señalando, por ejemplo, que se trató del único español en asistir repetidamente a congresos internacionales y que presentó comunicaciones en varios de ellos. También se mencionan sus colaboraciones en publicaciones internacionales, como su trabajo en el número fundacional de \textit{L’Enseignement Mathématique} \cite{FUR1}, que fue el delegado español en el CIEM a partir del congreso de Roma de 1908 o que fue miembro de la \textit{Commission Permanente du Répertoire Bibliographique des Sciences Mathématiques} \cite{AUS}.

La proyección internacional de García de Galdeano le llevó a formar parte de sociedades científicas y matemáticas extranjeras de su época. En una hoja de servicios conservada en el archivo histórico de la Universidad de Zaragoza (Signatura 18-A-2-6 (1), ver Figura 2) el propio García de Galdeano repasa las sociedades matemáticas a las que pertenecía hacia 1901.

\begin{figure}[h]
\label{Hoja}
\center
\includegraphics[width=7cm]{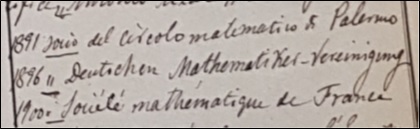}
\caption{Sociedades matemáticas de las que García de Galdeano era socio en 1901.}
\end{figure}

Sin embargo, no parece haberse prestado demasiada atención a esta pertenencia de García de Galdeano a varias de las sociedades matemáticas más importantes de su tiempo. En este trabajo tratamos de aportar luz sobre esta faceta de García de Galdeano presentando la información disponible encontrada en las revistas que actuaban como órganos de comunicación de dichas sociedades, así como en algunas otras fuentes documentales localizadas. 

Posiblemente el conocimiento y pertenencia a estas sociedades científicas extranjeras, influyera en su determinación para apoyar y colaborar en la fundación de la \textit{Sociedad Matemática Española} (SME). García de Galdeano estuvo involucrado en la organización de la Sección de Exactas del congreso fundacional de la \textit{Asociación para el Progreso de la Ciencias} celebrado en Zaragoza del 22 al 29 de septiembre de 1908.  Justamente en este congreso, a partir de la conferencia del general Benítez, toma fuerza y cuerpo la idea de fundar una sociedad matemática española similar a las existentes en otros países. Es más \cite[p.288]{GON}, en la sesión del día 25 de la Sección de Exactas presidida por V. Vera 

\begin{quote}
dio cuenta de la propuesta por la Comisión referente a la Asociación de Matemáticos. El Sr. Benítez abogó por que esta Sociedad funcionase independientemente de la Asociación general, haciendo además algunas observaciones los Sres. Silván, Jiménez Rueda y Presidente, y acordándose, finalmente, no prescindir de la antedicha Asociación, sino colocar aquélla bajo la protección de ésta, quedando encargado el Sr. Galdeano de formular la propuesta de la constitución de la Sociedad de Matemáticos que ha de presentarse al Comité Ejecutivo.
\end{quote}

Finalmente, la \textit{Sociedad Matemática Española} fue fundada en 1911, siendo su primer presidente D. José Echagaray y Eizaguirre\footnote{Para más detalles véanse, por ejemplo, los trabajos de Español \cite{ESP} o de Peralta \cite{PER}.}. Zoel García de Galdeano fue uno de los 359 socios fundadores de la SME. En la reunión de la Junta Directiva y de la General de la SME de 27 de octubre de 1913 se toma el siguiente acuerdo (\textit{Revista de la Sociedad Matemática Española}. Nº 22, noviembre de 1913, p. 59):

\begin{quote}
Por unanimidad se tomó el acuerdo de nombrar socio correspondiente al Sr. García de Galdeano, en atención á haber consagrado una gran parte de su vida al estudio y publicación de importantes obras y folletos de Matemáticas.
\end{quote}

De esta forma se empezaba a reconocer la titánica labor realizada durante varias décadas por García de Galdeano para modernizar las matemáticas españolas. El mayor reconocimiento por parte de la SME tuvo lugar a la muerte de José Echegaray en 1916. El 7 de diciembre de 1916, Zoel García de Galdeano fue nombrado segundo presidente de la SME, ocupando el cargo hasta su renuncia en 1920. En el discurso de aceptación del cargo, García de Galdeano menciona (\textit{Revista de la Sociedad Matemática Española}. Nº 51, enero de 1917, p. 50):

\begin{quote}
Entrando en materia, he de referirme para relacionarlo con la futura labor de nuestra Asociación matemática, a su desenvolvimiento progresivo, y hemos de contemplar el separatismo del pasado con el fusionismo del presente.
\end{quote}

Su presidencia fue más bien honorífica, aunque también tuvo aportaciones originales. Así, ante la refundación de la SME y la aprobación de unos nuevos estatutos en 1919, García de Galdeano, realiza el esfuerzo laboral y económico de editar el \textit{Suplemento a la Revista de la Sociedad Matemática Española}, de contenido crítico-bibliográfico, tal y como aparece recogido en el punto IX de los estatutos. Además, presidió las reuniones de la SME en Madrid en la visita que Jacques Hadamard hizo en marzo y abril de 1919 invitado por el Laboratorio de Matemáticos de la \textit{Junta de Ampliación de Estudios} dirigido por Julio Rey Pastor (Figura 3). 

\begin{figure}[h]
\center
\includegraphics[width=7cm]{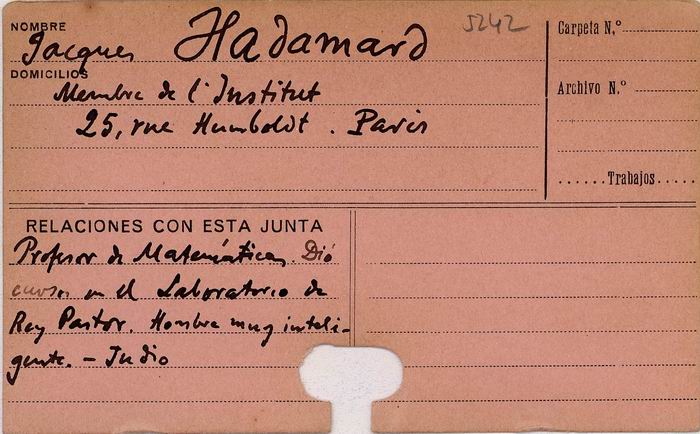}
\caption{Ficha de Jacques Hadamard. Fuente: Archivo de la Secretaría de la Junta para Ampliación de Estudios, Residencia de Estudiantes, Madrid.}
\end{figure}

\section{Zoel García de Galdeano en el \textit{Circolo Matematico di Palermo}}

El \textit{Circolo Matemático di Palermo} fue fundado el 2 de marzo de 1884 por el matemático de familia acomodada Giovanni Battista Guccia. Aunque el \textit{Circolo} no fue la primera sociedad matemática en fundarse, un cambio en los estatutos en febrero de1888 permitió la asociación de extranjeros en la sociedad. Tras ello, rápidamente se convirtió en la sociedad internacional más importante con la publicación matemática líder, los \textit{Rendiconti del Circolo Matemático di Palermo}. Figuras internacionales como Henri Poincaré, Gösta Mittag-Leffler, David Hilbert y Felix Klein pertenecieron a esta sociedad que en 1914 llegó a tener 924 socios frente a los 320 de la \textit{London Mathematical Society}, los 298 de la \textit{Societé Mathématique de France}; los 703 de la \textit{American Mathematical Society} o los 769 de la \textit{Deutsche Mathematiker-Vereinigung}.

\begin{figure}[h]
\center
\includegraphics[width=7cm]{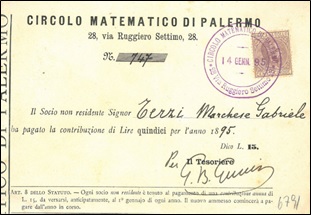}
\caption{Recibo de 1895 del Circolo de Palermo. Fuente: Colección privada de P.J. Miana.}
\end{figure}

Tras una votación secreta, Zoel García de Galdeano es nombrado socio no residente del \textit{Circolo} en la reunión del 24 de mayo de 1891\footnote{En 1891, el coste anual para los socios no residentes era de 15 Liras (aproximadamente el mismo coste que la suscripción anual de un diario) e incluía la subscripción gratuita a los \textit{Rendiconti} (Figura 4).} (\textit{Rendiconti del Circolo Matemático di Palermo}. Tomo 5, 1891, parte prima, p. 320). García de Galdeano fue presentado como socio por los secretarios  Guccia y Michele Luigi Albeggiani. En la misma sesión también se nombra socio no residente a Charles-Agne Laisant, presidente de la Societé Mathématique de France. El número de registro que le corresponde a García de Galdeano es el 177, siendo el primer español socio del \textit{Circolo}. En la siguiente reunión de 14 de junio, se registra la carta de agradecimiento del nuevo socio (\textit{Rendiconti del Circolo Matemático di Palermo}. Tomo 5, 1891, parte prima, p. 323). 

Durante los siguientes años las relaciones entre García de Galdeano y el \textit{Circolo} fueron principalmente bibliográficas. En la sección \textit{Biblioteca Matematica} de los \textit{Rendiconti} de los años 1892 a 1896 aparecen citadas varias obras de Galdeano en el apartado de publicaciones no periódicas (posiblemente enviadas por él mismo). Cada año, \textit{El Progreso Matemático} se incluye en el apartado de publicaciones periódicas con las que se realizaban intercambios, incluyendo los títulos y autores de los artículos publicados en el año anterior. Análogamente, en \textit{El Progreso Matemático} (Vol. 1, 20 de mayo de 1891, pp. 119-120) se incluye un pequeño artículo sobre el \textit{Circolo}. Este número está fechado cuatro días antes del nombramiento de García de Galdeano como socio no residente. Se menciona el contenido de los dos primeros fascículos del \textit{Rendiconti} de ese año 1891, y se describe la filosofía del \textit{Circolo} que permitía mediante intercambio de sus publicaciones acrecentar su \textit{Biblioteca Matematica}. Posiblemente fuera este uno de los antecedentes que animaron a García de Galdeano a fundar la primera revista matemática española.  Algunos de los artículos recogidos en los Rendiconti aparecen traducidos en \textit{El Progreso Matemático}. Por ejemplo, la necrológica de Enrico Betti firmada por Eugenio Beltrami fechada el 3 de enero de 1893 en Roma y que aparece en el tomo VI de los Rendiconti se encuentra traducida en \textit{El Progreso Matemático} (Vol. 25, 15 de enero de 1893, pp 30-31). En la edición de 1914 del \textit{Annuario Biografico del Circolo Matematico di Palermo} se incluye una pequeña biografía, posiblemente escrita por él mismo.

En 1914, al cumplirse el 30 de aniversario de la fundación del Circolo, se realiza un estudio detallado de la sociedad bajo el título de \textit{Nota Statistiche}. En él, Zoel García de Galdeano permanecía como socio no residente en activo junto con otros diez españoles: José Ríus y Casas (1904); José Gabriel Álvarez-Ude (1905); Lauro Clariana Ricart (1906); Esteban Terradas è Ylla (1909); Julio Rey Pastor (1911); Jorge Torner de la Fuente, Francisco Cebrián, Ricardo Cirera y Fernando Ruiz-Feduchy (1913); y Antonio Torroja y Miret (1914). Además, otros españoles habían sido socios y por motivos variados habían causado baja como, por ejemplo, Juan Jacobo Durán Loriga (1896); Angel Bozal Obejero (1903) y José Ruiz-Castizo (1907). Es destacable que los catedráticos de la Universidad Central Octavio de Toledo y Cecilio Jiménez Rueda ingresan en el Circolo en 1919, más de 25 años más tarde que García de Galdeano. (\textit{Annuario Biografico del Circolo Matematico di Palermo}, Palermo, 1928, pp. 115-126).

\section{Zoel García de Galdeano en la \textit{Deutsche Mathematiker-Vereinigung}}

La sociedad matemática alemana fue fundada en 1890. El órgano de comunicación oficial de la sociedad es el anuario que ésta viene publicando desde el año de su fundación \cite{SCH}.  
Tenemos constancia por una carta conservada en la biblioteca de la Universidad de Freiburg de que García de Galdeano ingresó, por invitación, en la \textit{Deutsche Mathematiker-Vereinigung} en junio de 1897. La carta en la que García de Galdeano acepta la invitación, cuyo destinatario es desconocido, está escrita en francés en papel con membrete de la Facultad de Ciencias de la Universidad de Zaragoza. A continuación transcribimos el texto de dicha carta:

\begin{quote}
Zaragoza 14 de Junio de 1897.\\
Monsieur:\\
Après une absence de deux mois à mon retour à Saragosse, j’ai eu le plaisir de trouver l’aimable invitation du \textit{Deutsche Mathematiker-Vereinigung} pour y apartenir, et m’empressant à vous envoyer mon acceptation, je vous prie de exprimer à la Commision le temoinage de ma plus profonde gratitude à invitation si flateuse pour moi.
Je vous adresse ci-joint un cheque de dix francs ; et je profitte cette ocassion de vous envoyer la temoignage de ma très haute estime\\
Zoel Gª de Galdeano
Coso 99 – 3º
\end{quote}

No tenemos constancia de quién cursó la invitación, ni bajo qué circunstancias. El presidente de la sociedad en 1897 era Felix Klein, con quién posteriormente García de Galdeano tendría relación al presidir éste el CIEM, pero sabemos por una carta de García de Galdeano que en 1899 ambos matemáticos no se conocían. Por otro lado, en agosto de 1897 tuvo lugar el congreso internacional de Zurich. Pese a que no asistió al mismo, en las actas de este congreso, se incluye un trabajo de 5 páginas de García de Galdeano titulado L’unification de concepts dans les mathématiques. Quizás exista alguna relación entre ambos hechos. 

El primer listado de socios que recoge el nombre de García de Galdeano apareció publicado en el vol. 6 del \textit{Jahresbericht der Deutschen Mathematiker-Vereinigung} (pp. 13-22). Entre los 374 socios que se enumeran, García de Galdeano es el único español. Hasta 1902 no se incorpora un nuevo socio español, José Rius y Casas, que lo fue hasta 1910. En el vol. 33 del \textit{Jahresbericht der Deutschen Mathematiker-Vereinigung} se da cuenta del fallecimiento de García de Galdeano y de su consiguiente baja como socio. En el listado de los 993 socios (entre personas e instituciones) a fecha 1 de mayo de 1926 que aparece en el vol. 35 del \textit{Jahresbericht der Deutschen Mathematiker-Vereinigung} sólo aparecen dos españoles: Julio Rey Pastor (socio desde 1913) y Tomás Rodríguez Bachiller (socio desde 1923).

En cualquier caso las interacciones oficiales entre García de Galdeano y la sociedad alemana debieron ser escasas. De hecho, las menciones a García de Galdeano en los anuarios de la sociedad son prácticamente inexistentes más allá de la aparición de algunas de sus obras en listados de bibliografía que en ningún caso fueron reseñadas en detalle.

\section{Zoel García de Galdeano en Francia. La AFAS y la Societé Mathématique de France}

En la sesión del 21 de noviembre de 1900, y bajo la presidencia de Henri Poincaré es elegido por unanimidad García de Galdeano como miembro de la Societé. Fue presentado por Charles-Agne Laisant y Émile Borel y continuó siendo socio hasta su muerte según podemos ver en la sección \textit{Vie de la Societé} (\textit{Bulletin de la Societé Mathématique de France}. Tome 29, 1901, p. 15 y Tome 51 (supplément spécial), 1923, p. 6). 
 
Sin embargo, García de Galdeano no fue el primer matemático español en pertenecer a la Societé. Se trató del Capitán Ignacio Beyens, que ingresó en 1887. Posteriormente lo hicieron Leonardo Torres Quevedo (1896) y Juan Jacobo Durán Loriga (1897). En 1902 ingresó el matemático gallego David Fernández Diéguez. 

Uno de los objetivos que se marcó la \textit{Societé Mathématique de France} en su fundación en 1885 fue la creación de un catálogo bibliográfico. El continuo crecimiento de las revistas matemáticas y de los artículos publicados animó a la \textit{Societé} a crear una herramienta útil para estudiantes y matemáticos profesionales. El \textit{Répertoire Bibliographique des Sciences Mathématiques} llegó a ser un desafío internacional. Dirigidos por Henri Poincaré durante casi 27 años (1885-1912), 50 matemáticos de 16 países revisaron más de 300 revistas matemáticas para sacar más de 20.000 referencias bibliográficas. En un principio se pensó publicar este catálogo en formato de libro, pero finalmente se decidió hacerlo en un formato de fichas (ver Figura 5) clasificadas por campos matemáticos\footnote{En la web http://sites.mathdoc.fr/RBSM/ puede consultarse los datos recogidos en estas fichas.}.

\begin{figure}[h]
\center
\includegraphics{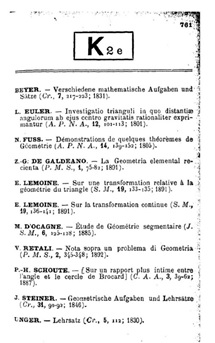}
\caption{Ficha 761 del Repertoire en la que se cita una obra de Zoel Garcia de Galdeano.
Fuente: http://sites.mathdoc.fr/RBSM/.}
\end{figure}

Aunque García de Galdeano no pertenecía a la \textit{Commission Permanente} nombrada en el \textit{Congrès Internatiònal de Bibliographie des Sciences Mathématiques} (16-19 de julio de 1889), entró a formar parte en ella el 1 de octubre de 1894. Por tanto, no parece adecuado atribuir a García de Galdeano el mérito de incluir el español (además del francés, inglés, alemán e italiano) entre los idiomas que la \textit{Commission} consideraría. Zoel García de Galdeano aparece en 37 entradas mientras que la revista \textit{El Progreso Matemático} lo hace en 68. También en el número 4 de \textit{El Progreso Matemático} (pp. 91-94) se presenta una traducción de las actas del congreso de 1889.

Uno de los primeros congresos internacionales a que acudió García de Galdeano fue a la 22ª Sesión de la \textit{L'Association française pour l'avancement des sciences} (AFAS) celebrado en Besançon del 3 al 10 de agosto de 1893. En esta reunión presentó la comunicación \textit{Note sur les institutions scientifiques et en particulier sur l’enseignement mathematique en Espagne}. 
La AFAS fue fundada en 1872 con el objetivo inicial de crear lazos de colaboraciones entre los investigadores principalmente franceses. García de Galdeano fue socio de la AFAS desde 1893 a 1897 y así se recoge en la lista anual de los \textit{Compte Rendu} de las sesiones 22 a 27 correspondientes a estos años. Otros españoles como Manuel Cano y de León, Gregorio Chil y Naranjo, Serafín Jordan, Ramón del Río ó Juan Vilanova y Piera eran ya miembros al ingresar García de Galdeano. 

García de Galdeano asistió  a la reunión de Saint-Étienne (Sesión 27, 5-12 de agosto de 1897) y a la de Paris (Sesión 30, del 2 al 9 de agosto de 1900) donde fue presidente de honor de la sección 16,  \textit{Pédagogie et Enseignement} e impartió la conferencia \textit{Quelques réflexions sur l'enseignement mathématique}.

En 1899 los matemáticos Charles-Agne Laisant (francés) y Henri Fehr (suizo) fundaron la revista \textit{L’Enseignement Mathématique} que con el tiempo se convirtió en la revista más importante de la enseñanza de las matemáticas. Los objetivos principales de la revista, según sus fundadores, eran la comunicación y la internacionalización. El Comité de Patronage, comité editorial, estaba formado por 20 matemáticos de 17 países diferentes.  El origen francés de esta iniciativa queda patente en los tres grandes nombres que forman la representación francesa: Henri Poincaré, Paul Appell y Émile Picard (los tres y Laisant fueron presidentes de la \textit{Societé Mathématique de France}). Otras figuras internacionales como Felix Klein, Moritz Cantor, o Gösta Mittag-Leffler acompañaron a Zoel García de Galdeano en este comité, siendo éste el único representante español. 

El primer artículo de la revista, después de la presentación de la revista por los directores, lleva por título \textit{Les mathemátiques en Espagne} y está firmado por Z. G. de Galdeano. Este artículo supone un recorrido preciso por las matemáticas españolas en el siglo XIX y refleja un amplio conocimiento de la situación matemática española de finales de ese siglo. 

La revista \textit{L’Enseignement Mathématique} promovió activamente el intercambio de ideas sobre la enseñanza de las matemáticas. En \cite{DES}, el estadounidense David Eugene Smith lanzaba la idea que una comisión internacional nombrada en un congreso internacional de matemáticos podía tratar adecuadamente la reforma de la enseñanza de las matemáticas puras y aplicadas, tanto en nivel de educación superior como universitario. Así en el Congreso Internacional de Matemáticos de Roma (1908) se tomó la resolución de nombrar la \textit{Commission Internationale de L’Enseignement Mathématique} (CIEM) siendo Felix Klein su presidente y Henri Fehr su secretario general.  La \textit{Commision} nombró delegados en cada país, siendo García de Galdeano el delegado español y presidente de la Subcomision española. Una de las primeras decisiones de la Commission fue adoptar la revista \textit{L’Enseignement Mathématique} como órgano oficial de comunicación.

Cuatro años más tarde, en el congreso de Cambridge, García de Galdeano fue cesado en esta representación con bastante polémica. La recién creada \textit{Sociedad Matemática Española} nombró en su puesto a Cecilio Jiménez Rueda. Actualmente esta comisión sigue funcionando con el nombre oficial de \textit{International Commission on Mathematical Instruction} (ICMI). 

El reconocimiento que los matemáticos franceses tuvieron por Zoel García de Galdeano se muestra en la carta que, fechada el 22 de junio de 1924 en Paris, Jacques Hadamard (Miembro Honorario de la SMF) envió a la Universidad de Zaragoza para manifestar sus condolencias por el reciente fallecimiento de García de Galdeano. Jacques Hadamard asistió a todos los congresos internacionales de matemáticos, coincidiendo con García de Galdeano. Como se ha comentado anteriormente, visitó Madrid en 1919 y posteriormente Barcelona en 1921. La carta mencionada, se reproduce en \cite[p. 359]{HOR4} aunque se encuentra ilocalizable actualmente. Dice así:

\begin{quote}
J’envoi à l’Université de Saragosse, à l’occasion de la perte qu’elle vient de faire en la personne du vénéré professeur Don Zoel de Galdeano, l’expression de ma profonde et douloreuse sympathie.\\ 
\begin{flushright}
J. Hadamard\\
25, rue Alexandre de Humboldt, Paris (14e)
\end{flushright}
\end{quote}

\section{Zoel García de Galdeano en la \textit{Mathematical Association of America}}

Desde 1894 ya existía en Estados Unidos una sociedad matemática, la \textit{American Mathematical Society}, que surgió a partir de la \textit{New York Mathematical Society} creada en 1888 \cite{ARC}. Pese a ello, una segunda sociedad, la \textit{Mathematical Association of America} fue fundada en 1915. Su órgano de comunicación fue desde sus inicios The \textit{American Mathematical Monthly}, revista que había sido creada muchos años antes, en 1894, y que todavía se edita actualmente \cite{BEN}. 
La primera interacción entre García de Galdeano y la MAA de la que tenemos noticia aparece recogida en el vol. 1, núm. 10 de octubre de 1894 del \textit{American Mathematical Monthly}. Allí se da cuenta de la recepción de los ocho primeros números del año 1894 de \textit{El Progreso Matemático} y se hace una breve, pero elogiosa, reseña de la revista firmada con las iniciales B. F. F. (probablemente de Benjamin Franklin Finker). En concreto, leemos (p.370): 

\begin{quote}
The editor of El Progresso [sic] Matemático has just sent us a copy of each of the eight issues for the year 1894. This journal is published monthly and is devoted to the solutions of problems in pure and applied mathematics. Many papers on interesting mathematical subjects also appear. The journal is well printed and many of the solutions are illustrated with beautiful black diagrams. The price of the journal to subscribers within the Postal Union is 11 fr.
\end{quote}

A lo largo de los números siguientes aparecieron algunas otras breves reseñas sobre \textit{El Progreso Matemático}. Por ejemplo, en el vol. 1, núm. 11 de noviembre de 1894 del \textit{American Mathematical Monthly}, con la firma J.M.C. (quizás John Marvin Colaw), leemos (p. 413): 

\begin{quote}
The October number of this excellent Monthly has several interesting papers, and many fine problems and solutions.
\end{quote}

El mismo texto y con la misma firma, pero referido al número de marzo de \textit{El Progreso Matemático}, aparece en el vol. 2, núm. 6 de junio de 1895 del \textit{Monthly}.
En el vol. 3, núm 3 de marzo de 1896 del \textit{Monthly}, B.F. Finker retoma los comentarios elogiosos con el siguiente texto (p. 94):

\begin{quote}
In this Journal are published problems which are proposed by the best mathematicians in the world. The solutions are illustrated with beautiful diagrams.
\end{quote}

Finalmente, en el vol. 3, núm 4 de abril de 1896 del \textit{Monthly} (p. 126) se acusa recibo de la recepción del Tomo V de \textit{El Progreso Matemático}, correspondiente al año 1895.
Las menciones a \textit{El Progreso Matemático} no se limitaron a breves descripciones de su contenido. Por ejemplo, en el vol. 2, núms. 7-8 de julio y agosto de 1895 del \textit{Monthly} aparece (p. 247) una nota comunicada por D. E. Smith en la que, bajo el título “The \textit{American Mathematical Monthly} in Spain”, se informa de que en el número de marzo de \textit{El Progreso Matemático} se menciona una serie de artículos sobre geometría no-euclídea que Bruce Halsted había venido publicando en el \textit{Monthly}.
Pero \textit{El Progreso Matemático} no fue la única publicación de García de Galdeano que mereció espacio en las páginas del Monthly. En el vol. 24, núm. 7 de septiembre de 1917 D.A. Rothrock escribe (p. 350):

\begin{quote}
The Spanish Mathematical Society has published a Supplement written by the president of the Society, Dr. Zoel García de Galdeano, upon the development of mathematics under the title “Exposición sumaria de la matemática según un nuevo método”.
\end{quote}

En el vol. 26, núm. 3 de marzo de 1919, esta vez sin firma leemos:

\begin{quote}
The first number of a new scientific periodical, El Progreso Científico, has appeared recently at Saragossa under the direction of Professor Z. G. de Galdeano. It is to be a semi-annual review devoted to mathematics, physics and chemistry, and containing papers dealing with fundamental questions, with criticism and with scientific methodology, as well as with bibliographies, and matters pertaining to the teaching of science.
\end{quote}

En otra de las reseñas que apareció en el vol. 26, núm. 7 de septiembre de 1919 del \textit{Monthly} (p. 302), no sólo se menciona la publicación en cuestión, sino que re hace un breve repaso a las revistas editadas por García de Galdeano:

\begin{quote}
We have also received Tome 1, no. 1 of Professor Z. G. de Galdeano new periodical Supplemento a la Revista Matemática Hispano-Americana, boletín de crítica, pedagogía, historia y bibliografía (32 pages). Professor Galdeano has also written or edited four other periodicals: (1) \textit{El Progreso Matemático}, 7 vols., 1891-1895, 1899-1900; (2) Boletín de crítica, enseñanza y bibliografía, 2 nos., 1907-1908; (3) Supplemento a la revista de la Sociedad Matemática Española, 3 nos., 1917; (4) El Progreso scientifico, revista semestral, 1 no., July 1918.
\end{quote}

Además de sus publicaciones, la propia figura de García de Galdeano fue mencionada en dos ocasiones a raíz de sendos congresos internacionales, el de Zurich de 1897 y el de París de 1900. Respecto al de Zurich, es interesante señalar que en la breve crónica a cargo de G. B. Halsted que apareció en el vol. 4, núms. 8-9 de agosto y septiembre de 1897 de Monthly el único ponente que se menciona explícitamente al repasar la sección segunda del congreso (dedicada al análisis y teoría de funciones) es García de Galdeano. En concreto (p. 229) se dedican a García de Galdeano las siguientes significativas palabras:

\begin{quote}
The second section contained a title from Z. De Galdeano, whose heroic efforts gave Spain a Journal of Mathematics, now unfortunately dead in the decadence of that beautiful priest-ridden land.
\end{quote}

En el caso del Congreso Internacional de París de 1900, la crónica aparecida en el vol. 7, núms. 8-9 de agosto y septiembre de 1900 corrió nuevamente a cargo de Halsted. Al inicio de la crónica se menciona al presidente del congreso (Poincaré), a algunos vicepresidentes y a los presidentes de las distintas secciones. De los delegados oficiales sólo tres son mencionados explícitamente por Halsted, uno de ellos es García de Galdeano.  

\begin{figure}[h]
\center
\includegraphics[width=7cm]{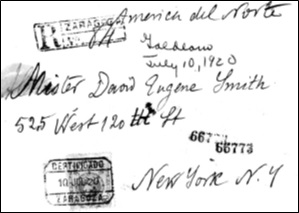}
\caption{Sobre de la carta de García de Galdeano a D.E. Smith. Fuente: David Eugene Smith Professional Papers, Rare Book \& Manuscript Library, Columbia University in the City of New York.}
\end{figure}

Sin embargo, pese a esta presencia relativamente abundante en la publicación oficial de la MAA, no es hasta 1920 que García de Galdeano pasa a ser socio de dicha asociación. La fecha exacta en la que fue aceptado como miembro fue el 6 de septiembre de1920. D hecho, en el vol. 27, núm. 11 de noviembre de 1920 del \textit{Monthly} aparece (pp. 389-390) un listado de las 70 personas aceptadas como miembros en la reunión del consejo de la asociación que tuvo lugar en esa fecha. En ese listado aparece un segundo español: Luis Octavio de Toledo. Respecto a este hecho, se puede afinar todavía más puesto que en la \textit{Rare Book \& Manuscript Library} de la Universidad de Columbia se conserva entre los papeles de D. E. Smith la carta en la que García de Galdeano solicita su admisión como socio. Dicha carta, escrita en castellano, está fechada en Zaragoza a 10 de julio de 1920 (ver Figura 6) y tiene sello de entrada en los EE.UU. el 27 de julio del mismo año. A continuación transcribimos el breve texto de la carta:

\begin{quote}
Sr. David Eugene Smith\\
Muy distinguido colega. Recibí su atenta carta y atendiendo a sus indicaciones le remito un cheque de 5 Dollars importe de mi inscripción como socio durante el primer año. Próximamente tendré el gusto de remitirle algunas de mis últimas publicaciones y le agradeceré que las cite en la revista de la sociedad. Quedando de V afectísimo amigo y colega.
Zoel G. de Galdeano\\
Zaragoza, Cervantes 5, pral. dcha.\\
10 julio 1920
\end{quote}

Aparentemente García de Galdeano permaneció como socio de la MAA hasta su muerte. Su nombre aparece en el listado de miembros publicado en el vol. 29, núm. 6 de junio y julio de 1922 del \textit{Monthly}, pero ya no en el listado de miembros publicado en el vol. 31, núm. 9 de noviembre de 1924 del \textit{Monthly}; ya después de su muerte. Desafortunadamente no hay registros intermedios de socios. Es importante destacar que, según el listado del año 1922 de los 1476 miembros de la MAA sólo 54 eran extranjeros (aparte de los canadienses) y, como hemos dicho, 2 españoles. Con la muerte de García de Galdeano, Octavio de Toledo fue el único socio español hasta 1934, año de su fallecimiento. Habrá que esperar hasta 1947 para volver a encontrar socios españoles en los listados de miembros de la MAA.

\section{Conclusiones}

La vocación internacional de Zoel García de Galdeano, que se ha puesto de manifiesto en múltiples ocasiones, queda también reflejada en su pertenencia a las sociedades profesionales de matemáticas más importantes de su época. Aunque sólo en algunos casos era el primer español en ingresar en esa sociedad, sí fue el primero en pertenecer simultánea e ininterrumpidamente a todas ellas.

Resulta fácil imaginar que la pertenencia a todas estas sociedades resultó crucial para que García de Galdeano tuviera noticia de los avances que se llevaban a cabo en las matemáticas de su tiempo, tanto en lo referente a investigación como en lo tocante a la enseñanza. Además, el acceso a publicaciones periódicas le permitiría conocer un buen número de obras a las que, de otro modo, no habría tenido acceso. El conocimiento adquirido con estas lecturas y con los intercambios realizados tanto a través de \textit{El Progreso Matemático}, como a título personal, explica en parte la gran labor introductoria divulgadora que García de Galdeano llevó a cabo en España durante casi medio siglo.

Finalmente, la experiencia ganada con los años de pertenencia a alguna de estas asociaciones debió de resultarle muy útil en los inicios de la \textit{Sociedad Matemática Española} y, en particular, durante su presidencia de la misma.

\section*{Agradecimientos}
Esta investigación ha sido parcialmente financiada por el proyecto FCT-16-1216: ``Año García de Galdeano, pasión por las matemáticas''. Los autores agradecen al profesor Luis Español diversos comentarios, detalles y sugerencias que han ayudado a la realización de este artículo, así como a Alfredo Valverde por remitirnos una copia de la ficha de Hadamard en la Residencia de Estudiantes.
\end{document}